\documentclass[11pt, twoside]{article}
\usepackage{latexsym}
\usepackage{amsmath}
\usepackage{amssymb}
\usepackage[all]{xy}
\usepackage{amsfonts}
\usepackage{verbatim}
\usepackage{amsthm}
\usepackage{mathrsfs}
\usepackage{epsfig}
\usepackage{xy}
\usepackage{array}
\usepackage{stmaryrd}
\usepackage{graphicx,color}
\usepackage{xcolor}
\usepackage[colorlinks=true,linkcolor=blue,citecolor=blue]{hyperref}
\usepackage{tikz}
\usetikzlibrary{arrows,calc}
\usepackage{etex}
\usepackage{mathdots}
\usepackage{float}
\usepackage{graphics}
\usepackage{pdflscape}
\usepackage{bm}
\usepackage{anysize,hyperref}
\input xypic
\xyoption{all}
\newcommand{\mr}{\hbox{\boldmath$\cdot$}}
\usepackage[perpage,symbol]{footmisc}
\usepackage{setspace}
\def\A{\mathscr{A}}

\def\C{\mathscr{C}}
\def\E{\mathbb{E}}
\def\F{\mathbb{F}}
\def\s{\mathfrak{s}}
\def\id{\mathrm{id}}
\def\op{^\mathrm{op}}
\def\Ab{\mathit{Ab}}
\def\del{\delta}
\def\dr{\ar@{->}[r]}

\def\X{\mathscr{X}}

\newcommand{\CC}{{\bf{C}}^{n+2}_{\C}}

\newcommand{\ov}{\overset}
\newcommand{\lra}{\longrightarrow}
\newcommand{\co}{\colon}
\newcommand{\uas}{^{\ast}}            
\newcommand{\sas}{_{\ast}}
\newcommand{\Xd}{\langle X^{\mr},\del\rangle}  
\newcommand{\Yr}{\langle Y^{\mr},\rho\rangle}  

\SelectTips{cm}{10}

\begin{document}
\title{\Large{\bf $\bm{n}$-extension closed subcategories of $\bm{(n+2)}$-angulated categories$^\bigstar$\footnotetext{\hspace{-1em}$^\bigstar$This work was supported by the National Natural Science Foundation of China (Grant No. 11901190) and the Scientific Research Fund of Hunan Provincial Education Department (Grant No. 19B239).}}}
\medskip
\author{Panyue Zhou}

\date{}

\maketitle
\def\blue{\color{blue}}
\def\red{\color{red}}

\newtheorem{theorem}{Theorem}[section]
\newtheorem{lemma}[theorem]{Lemma}
\newtheorem{corollary}[theorem]{Corollary}
\newtheorem{proposition}[theorem]{Proposition}
\newtheorem{conjecture}{Conjecture}
\theoremstyle{definition}
\newtheorem{definition}[theorem]{Definition}
\newtheorem{question}[theorem]{Question}
\newtheorem{remark}[theorem]{Remark}
\newtheorem{remark*}[]{Remark}
\newtheorem{example}[theorem]{Example}
\newtheorem{example*}[]{Example}
\newtheorem{condition}[theorem]{Condition}
\newtheorem{condition*}[]{Condition}
\newtheorem{construction}[theorem]{Construction}
\newtheorem{construction*}[]{Construction}

\newtheorem{assumption}[theorem]{Assumption}
\newtheorem{assumption*}[]{Assumption}

\baselineskip=17pt
\parindent=0.5cm
\vspace{-6mm}

\begin{abstract}
\baselineskip=16pt
Let $\C$ be a Krull-Schmidt $(n+2)$-angulated category and $\A$ be an $n$-extension closed subcategory of $\C$.
Then $\A$ has the structure of an $n$-exangulated category in the sense of
Herschend--Liu--Nakaoka. This construction gives $n$-exangulated categories which are not
$n$-exact categories in the sense of Jasso nor $(n+2)$-angulated categories  in the sense of
 Geiss--Keller--Oppermann in general.
As an application, our result can lead to a recent main result of Klapproth.\\[0.2cm]
\end{abstract}

\pagestyle{myheadings}
\markboth{\rightline {\scriptsize   Panyue Zhou}}
         {\leftline{\scriptsize $n$-extension closed subcategories of $(n+2)$-angulated categories}}

\section{Introduction}
In \cite{GKO}, Geiss, Keller and Oppermann introduced $(n+2)$-angulated categories. These
are a ``higher dimensional" analogue of triangulated categories,
 in the sense that triangles are replaced by $(n+2)$-angles, that is, morphism sequences of length $(n+2)$. Thus a $3$-angulated category is precisely a triangulated category.
An important source of examples of $(n+2)$-angulated categories are certain cluster tilting subcategories of triangulated categories.
Jasso \cite{Ja} introduced $n$-exact categories as
higher analogs of exact categories. Moreover, he also proved that any $n$-cluster-tilting subcategory of an
exact category is an $n$-exact category.

 Let  $(\C,\Sigma,\Theta)$ be a  Krull-Schmidt $(n+2)$-angulated
category and $\A$ be an $n$-extension closed subcategory of $\C$.
An $\A$-conflation is a complex
$$A_0\xrightarrow{f_0}A_1\xrightarrow{f_1}A_2\xrightarrow{f_2}\cdots\xrightarrow{f_{n-1}}A_n\xrightarrow{f_n}A_{n+1}$$
with $A_0,A_1,\cdots, A_{n+1}\in\A $ for which there exists a morphism $f_{n+1}\colon A_{n+1}\to \Sigma A_0$ such that
$$A_0\xrightarrow{f_0}A_1\xrightarrow{f_1}A_2\xrightarrow{f_2}\cdots\xrightarrow{f_{n-1}}A_n\xrightarrow{f_n}A_{n+1}\xrightarrow{f_{n+1}}\Sigma A_0$$
is an $(n+2)$-angle in $(\C,\Sigma,\Theta)$. We denote by $\mathscr{E}_{\A}$ the class of all
$\A$-conflations.

Klapproth proved the following result.

\begin{theorem}{\rm\cite[Theorem 3.2]{K}}\label{main}
Let  $(\C,\Sigma,\Theta)$ be a  Krull-Schmidt $(n+2)$-angulated
category and $\A$ be an $n$-extension closed subcategory of $\C$.
If $\C(\Sigma\A,\A)=0$, then $(\A,\mathscr{E}_{\A})$ is an $n$-exact category.
\end{theorem}
In particular, when $n=1$, it is the main theorem of \cite{D} and it was also rediscovered in \cite[Proposition 2.5]{J}.

Recently, Herschend, Liu and Nakaoka \cite{HLN} introduced the
notion of $n$-exangulated categories for any positive integer $n$. It is not only a higher
dimensional analogue of extriangulated categories defined by Nakaoka and Palu \cite{NP},
but also gives a common generalization of $n$-exact categories in the sense of
Jasso \cite{Ja} and $(n+2)$-angulated categories in the sense of Geiss-Keller-Oppermann \cite{GKO}.

Nakaoka and Palu \cite{NP} proved extension closed subcategories of triangulated categories are
extriangulated categories. This construction gives extriangulated categories which are not
exact and triangulated.
Based on this idea, we prove the following conclusion.
\begin{theorem}\label{main2}{\rm (see Theorem \ref{main1} for details)}
Let  $(\C,\Sigma,\Theta)$ be a  Krull-Schmidt $(n+2)$-angulated
category and $\A$ be an $n$-extension closed subcategory of $\C$.
Then $\A$ has the structure of an $n$-exangulated category, induced from that of $\C$.
\end{theorem}

From this theorem, we know that any $(n+2)$-angulated category can be viewed as
an $n$-exangulated category. This generalizes the Proposition 4.5 in \cite{HLN} and
also is a higher couterpart of Nakaoka-Palu's result.

Herschend, Liu and Nakaoka \cite[Proposition 4.37]{HLN}  gave a description of when an $n$-exangulated category can become $n$-exact category.
Based on this fact and Theorem \ref{main}, we give a new proof of Theorem \ref{main}.
Our proof method is to avoid proving that $\mathscr{E}_{\A}$ is closed under weak isomorphisms, which is very complicated and difficult.

This article is organized as follows. In Section 2, we review some elementary definitions and facts on $n$-exangulated categories.
In Section 3, we prove our main result in this article.

\section{Preliminaries}
In this section, let $\C$ be an additive category and $n$ be a positive integer. Suppose that $\C$ is equipped with an additive bifunctor $\E\colon\C\op\times\C\to{\rm Ab}$, where ${\rm Ab}$ is the category of abelian groups. Next we briefly recall some definitions and basic properties of $n$-exangulated categories from \cite{HLN}. We omit some
details here, but the reader can find them in \cite{HLN}.

{ For any pair of objects $A,C\in\C$, an element $\del\in\E(C,A)$ is called an {\it $\E$-extension} or simply an {\it extension}. We also write such $\del$ as ${}_A\del_C$ when we indicate $A$ and $C$. The zero element ${}_A0_C=0\in\E(C,A)$ is called the {\it split $\E$-extension}. For any pair of $\E$-extensions ${}_A\del_C$ and ${}_{A'}\del{'}_{C'}$, let $\delta\oplus \delta'\in\mathbb{E}(C\oplus C', A\oplus A')$ be the
element corresponding to $(\delta,0,0,{\delta}{'})$ through the natural isomorphism $\mathbb{E}(C\oplus C', A\oplus A')\simeq\mathbb{E}(C, A)\oplus\mathbb{E}(C, A')
\oplus\mathbb{E}(C', A)\oplus\mathbb{E}(C', A')$.

For any $a\in\C(A,A')$ and $c\in\C(C',C)$,  $\E(C,a)(\del)\in\E(C,A')\ \ \text{and}\ \ \E(c,A)(\del)\in\E(C',A)$ are simply denoted by $a_{\ast}\del$ and $c^{\ast}\del$, respectively.

Let ${}_A\del_C$ and ${}_{A'}\del{'}_{C'}$ be any pair of $\E$-extensions. A {\it morphism} $(a,c)\colon\del\to{\delta}{'}$ of extensions is a pair of morphisms $a\in\C(A,A')$ and $c\in\C(C,C')$ in $\C$, satisfying the equality
$a_{\ast}\del=c^{\ast}{\delta}{'}$.}
Then the functoriality of $\E$ implies $\E(c,a)=a_{\ast}(c^{\ast}\del)=c^{\ast}(a_{\ast}\del)$.

\begin{definition}\cite[Definition 2.7]{HLN}
Let $\bf{C}_{\C}$ be the category of complexes in $\C$. As its full subcategory, define $\CC$ to be the category of complexes in $\C$ whose components are zero in the degrees outside of $\{0,1,\ldots,n+1\}$. Namely, an object in $\CC$ is a complex $X^{\mr}=\{X_i,d^X_i\}$ of the form
\[ X_0\xrightarrow{d^X_0}X_1\xrightarrow{d^X_1}\cdots\xrightarrow{d^X_{n-1}}X_n\xrightarrow{d^X_n}X_{n+1}. \]
We write a morphism $f^{\mr}\co X^{\mr}\to Y^{\mr}$ simply $f^{\mr}=(f^0,f^1,\ldots,f^{n+1})$, only indicating the terms of degrees $0,\ldots,n+1$.
\end{definition}

\begin{definition}\cite[Definition 2.23]{HLN}
Let $\s$ be an exact realization of $\E$.
\begin{enumerate}
\item[\rm (1)] An $n$-exangle $\Xd$ is called an $\s$-{\it distinguished} $n$-exangle if it satisfies $\s(\del)=[X^{\mr}]$. We often simply say {\it distinguished $n$-exangle} when $\s$ is clear from the context.
\item[\rm (2)]  An object $X^{\mr}\in\CC$ is called an {\it $\s$-conflation} or simply a {\it conflation} if it realizes some extension $\del\in\E(X_{n+1},X_0)$.
\item[\rm (3)]  A morphism $f$ in $\C$ is called an {\it $\s$-inflation} or simply an {\it inflation} if it admits some conflation $X^{\mr}\in\CC$ satisfying $d_X^0=f$.
\item[\rm (4)]  A morphism $g$ in $\C$ is called an {\it $\s$-deflation} or simply a {\it deflation} if it admits some conflation $X^{\mr}\in\CC$ satisfying $d_X^n=g$.
\end{enumerate}
\end{definition}

\begin{definition}\cite[Definition 2.32]{HLN}
An {\it $n$-exangulated category} is a triplet $(\C,\E,\s)$ of additive category $\C$, additive bifunctor $\E\co\C\op\times\C\to\Ab$, and its exact realization $\s$, satisfying the following conditions.
\begin{itemize}
\item[{\rm (EA1)}] Let $A\ov{f}{\lra}B\ov{g}{\lra}C$ be any sequence of morphisms in $\C$. If both $f$ and $g$ are inflations, then so is $g\circ f$. Dually, if $f$ and $g$ are deflations, then so is $g\circ f$.

\item[{\rm (EA2)}] For $\rho\in\E(D,A)$ and $c\in\C(C,D)$, let ${}_A\langle X^{\mr},c\uas\rho\rangle_C$ and ${}_A\Yr_D$ be distinguished $n$-exangles. Then $(\id_A,c)$ has a {\it good lift} $f^{\mr}$, in the sense that its mapping cone gives a distinguished $n$-exangle $\langle M^{\mr}_f,(d^X_0)\sas\rho\rangle$.
 \item[{\rm (EA2$\op$)}] Dual of {\rm (EA2)}.
\end{itemize}
Note that the case $n=1$, a triplet $(\C,\E,\s)$ is a  $1$-exangulated category if and only if it is an extriangulated category, see \cite[Proposition 4.3]{HLN}.
\end{definition}

\begin{example}
From \cite[Proposition 4.34]{HLN} and \cite[Proposition 4.5]{HLN},  we know that $n$-exact categories and $(n+2)$-angulated categories are $n$-exangulated categories.
There are some other examples of $n$-exangulated categories
 which are neither $n$-exact nor $(n+2)$-angulated, see \cite{HLN,LZ,HZZ}.
\end{example}

Let $(\C,\E,\s)$ be an $n$-exangulated category and  $\F\subseteq\E$ be an additive subfunctor {\rm (see \cite[Definition 3.7]{HLN})}. For a realization $\s$ of $\E$, define
$\s\hspace{-1.3mm}\mid_{\F}$ to be the restriction of $\s$ onto $\F$. Namely, it is defined by $\s\hspace{-1.3mm}\mid_{\F}(\del)=\s(\del)$ for any $\F$-extension $\del$.

\begin{lemma}\label{lem1}{\rm\cite[Proposition 3.16]{HLN}}
Let $(\C,\E,\s)$ be an $n$-exangulated category. For any additive subfunctor $\F\subseteq\E$, the following statements are equivalent.

{\rm (1)}~ $(\C,\F,\s\hspace{-1.2mm}\mid_{\F})$ is an $n$-exangulated category.

{\rm (2)}~ $\s\hspace{-1.2mm}\mid_{\F}$-inflations are closed under composition.

{\rm (3)}~ $\s\hspace{-1.2mm}\mid_{\F}$-deflations are closed under composition.

\end{lemma}

\section{Main result}
In this section, when we say that $\A$ is a subcategory of an additive category $\C$, we always assume that $\C$ is full, and closed under isomorphisms, direct sums and direct summands.

We first recall the notion of $n$-extension closed from \cite{L1}.

\begin{definition}\cite[Definition 3.6]{L1}
Let $(\C,\Sigma,\Theta)$ be an $(n+2)$-angulated category. A subcategory $\A$ of $\C$ is called \emph{$n$-extension closed} if
 for each morphism $f_{n+1}\colon A_{n+1}\to \Sigma^n A_0$ with $A_0,A_{n+1}\in\A$,
there exists an $(n+2)$-angle
$$A_0\xrightarrow{f_0}A_1\xrightarrow{f_1}A_2\xrightarrow{f_2}\cdots\xrightarrow{f_{n-1}}A_n\xrightarrow{f_n}A_{n+1}\xrightarrow{f_{n+1}}\Sigma A_0$$
with terms $A_1,A_2,\cdots,A_n\in\A$.
\end{definition}

The following lemma can be found in \cite[Lemma 2.5]{L2}.

\begin{lemma}\label{y1}{\rm\cite[Lemma 2.5]{L2}}
Let $(\C,\Sigma,\Theta)$ be a Krull-Schmidt $(n+2)$-angulated category and
$$A_{\bullet}:~A_0\xrightarrow{\binom{f_0}{g_0}}A_1\oplus B_1\xrightarrow{(f_1,\hspace{0.8mm} g_1)}A_2\xrightarrow{f_2}A_3\xrightarrow{f_3}\cdots\xrightarrow{f_{n-1}}A_n\xrightarrow{f_n}A_{n+1}\xrightarrow{f_{n+1}}\Sigma A_0.$$
is an $(n+2)$-angle in $\C$.
If $g_0=0$, then $A_{\bullet}\simeq A'_{\bullet}\oplus B_{\bullet}$, where
$$A'_{\bullet}:~A_0\xrightarrow{f_0}A_1\xrightarrow{f_{11}}A'_2\xrightarrow{f_{21}}A_3\xrightarrow{f_3}\cdots\xrightarrow{f_{n-1}}A_n\xrightarrow{f_n}A_{n+1}\xrightarrow{f_{n+1}}\Sigma A_0$$
and
$$B_{\bullet}:~0\xrightarrow{}B_1\xrightarrow{~1~}B_1\xrightarrow{}0\xrightarrow{}\cdots\xrightarrow{}0\xrightarrow{}0\xrightarrow{}0$$
is two an $(n+2)$-angle.
\end{lemma}

Let $(\C,\Sigma,\Theta)$ be an $(n+2)$-angulated category. Since $\Sigma\colon\C\xrightarrow{~\simeq~}\C$ is an  automorphism, then
$\Sigma$ gives an additive bifunctor
$$\E_{\Sigma}=\C(-,\Sigma-)\colon \C^{\rm op}\times \C\to {\rm Ab},$$
defined by the following.
\begin{itemize}
\item[\rm (i)] For any $A,C\in\C$, $\E_{\Sigma}(C, A)=\C(C,\Sigma A)$;

\item[\rm (ii)] For any $a\in\C(A,A')$ and $c\in\C(C',C)$, the map $\E_{\Sigma}(c, a)\colon\C(C, \Sigma A)\to \C(C', \Sigma A)$
sends $\delta\in\C(C, \Sigma A)$ to $c^{\ast}a_{\ast}\delta=(\Sigma a)\circ\delta\circ c$.
\end{itemize}

For each $\delta\in\E_{\Sigma}(C, A)$, we complete
it into an $(n+2)$-angle
$$A\xrightarrow{f_0}X_1\xrightarrow{f_1}X_2\xrightarrow{f_2}\cdots\xrightarrow{f_{n-1}}X_n\xrightarrow{f_n}C\xrightarrow{\delta}\Sigma A_0$$
Define $\s_{\Theta}(\delta)=[X^{\mr}]$ by using $X^{\mr}\in{\mathbf{C}^{n+2}_{(A,\hspace{0.8mm}C)}}$ given by
$$A\xrightarrow{f_0}X_1\xrightarrow{f_1}X_2\xrightarrow{f_2}\cdots\xrightarrow{f_{n-1}}X_n\xrightarrow{f_n}C$$

The following result shows that any $(n+2)$-angulated category can be viewed as
an $n$-exangulated category.
\begin{theorem}{\rm \cite[Proposition 4.5]{HLN}}
With the above definition, $(\C,\E_{\Sigma},\s_{\Theta})$ is an $n$-exangulated category.
\end{theorem}

Let $(\C,\Sigma,\Theta)$ be an $(n+2)$-angulated category and
$\A$ be an $n$-extension closed subcategory of $\C$.
We define $\E_{\A}$ to be the restriction of
$\E_{\Sigma}$ onto $\A^{\rm op}\times\A$, and define $\s_{\A}$ by restricting $\s_{\Theta}$.

The following construction gives $n$-exangulated categories which are not
$n$-exact  nor $(n+2)$-angulated in general.

\begin{theorem}\label{main1}
Let $(\C,\Sigma,\Theta)$ be a Krull-Schmidt $(n+2)$-angulated category and
$\A$ be an $n$-extension closed subcategory of $\C$. Then $(\A,\E_{\A},\s_{\A})$
is an $n$-exangulated category.
\end{theorem}

\proof It is straightforward to verify that $\E_{\A}$ is an additive subfunctor of $\E_{\Sigma}$.
By \cite[Lemma 3.8]{K}, we know that that $\s_{\A}$-inflations are closed under composition.
By Lemma \ref{lem1}, we have that $(\A,\E_{\A},\s_{\A})$
is an $n$-exangulated category. \qed

\begin{proposition}{\rm \cite[Proposition 4.37]{HLN}}\label{prop2}
Let $(\C, \E, \s)$  be an $n$-exangulated category. Assume that any $\s$-inflation is
monomorphic, and any $\s$-deflation is epimorphic in $\C$. Note that this is equivalent to
assuming that any $\s$-conflation is $n$-exact sequence. We denote the class of all $\s$-conflations by $\X$. If $(\C, \E, \s)$  satisfies the following conditions {\rm (a)} and {\rm (b)} for any pair of morphisms
$A\xrightarrow{~a~}B\xrightarrow{~b~}C$ in $\C$, then $(\C, \X)$  becomes an
n-exact category in the sense of Jasso {\rm \cite[Definition 4.2]{Ja}}.

{\rm (a)} If $b\circ a$ is an $\s$-inflation, then so is $a$.

{\rm (b)} If $b\circ a$  is an $\s$-deflation, then so is $b$.

\end{proposition}

In Theorem \ref{main1}, we know that $(\A,\E_{\A},\s_{\A})$
is an $n$-exangulated category.
We assume that any $\s_{\A}$-inflation is monomorphic, and any $\s_{\A}$-deflation is epimorphic in $\A$. Note that this is equivalent to assuming that any $\s$-conflation is $n$-exact sequence. Thus we denote the class of all $\s_{\A}$-conflations by $\mathscr{E}_{\A}$.

Now we give a new proof of Klapproth's result.

\begin{theorem}{\rm\cite[Theorem 3.2]{K}}\label{cor1}
Let $(\C,\Sigma,\Theta)$ be a Krull-Schmidt $(n+2)$-angulated category and
$\A$ be an $n$-extension closed subcategory of $\C$.
If $\C(\Sigma\A,\A)=0$, then $(\A,\mathscr{E}_{\A})$
is an $n$-exact category.
\end{theorem}

\proof By Theorem \ref{main1}, we know that $(\A,\E_{\A},\s_{\A})$
is an $n$-exangulated category.

We first claim that any $\s_{\A}$-inflation is monomorphic $\A$.

Assume that $f\colon A\to B$ is an $\s_{\A}$-inflation in $\A$.
By definition  there exists an $(n+2)$-angle
$$A\xrightarrow{~f~}B\xrightarrow{g_1}C_2\xrightarrow{g_2}\cdots\xrightarrow{g_{n-1}}C_n\xrightarrow{g_n}C_{n+1}\xrightarrow{g_{n+1}}\Sigma A$$
with terms $C_2,C_3,\cdots,C_n\in\A$.

Now we prove that $f$ is monomorphism in $\A$.
Let $h\colon C\to A$ be a morphism in $\A$ such that $fh=0$.
Apply the functor $\C(C,-)$ to the above  $(n+2)$-angle, we have the following exact sequence:
$$\C(C,\Sigma^{-1}C_{n+1})\xrightarrow{~}\C(C,A)\xrightarrow{\C(C,\hspace{0.8mm}f)}\C(C,B)$$
where the leftmost term vanishes since $\C(C,\Sigma^{-1}C_{n+1})\simeq \C(\Sigma C,C_{n+1})=0$.
This shows that $\C(C,f)$ is a monomorphism which implies $h=0$ since $fh=0$.
This shows that $f$ is monomorphism in $\A$.

In a similar way, one can show that any $\s_{\A}$-deflation is epimorphic in $\A$.

For any pair of morphisms
$A\xrightarrow{~a~}B\xrightarrow{~b~}C$ in $\A$,
we claim that if $b\circ a$ is an $\s_{\A}$-inflation, then so is $a$.
Indeed, since $h_0:=b\circ a$ is an $\s_{\A}$-inflation, by definition  there exists an $(n+2)$-angle
$$A\xrightarrow{~h_0~}C\xrightarrow{h_1}X_2\xrightarrow{h_2}\cdots\xrightarrow{h_{n-1}}X_n\xrightarrow{h_n}X_{n+1}\xrightarrow{h_{n+1}}\Sigma A$$
with terms $X_2,X_3,\cdots,X_n\in\A$.

For the morphism $X_{n+1}\xrightarrow{d_{n+1}:=\Sigma a\circ h_{n+1}}\Sigma B$ with $B,C\in\A$,
since $\A$ is $n$-extension closed, then there exists an $(n+2)$-angle
$$B\xrightarrow{d_0}Y_1\xrightarrow{d_1}Y_2\xrightarrow{d_2}
\cdots\xrightarrow{d_{n-1}}Y_{n}\xrightarrow{d_{n}}X_{n+1}\xrightarrow{d_{n+1}}\Sigma B$$
with $Y_1,Y_2,\cdots,Y_{n}\in\A$.
Consider the following commutative diagram of $(n+2)$-angles
$$\xymatrix{
\Sigma^{-1}X_{n+1}\ar[r]^{\quad h_{-1}}\ar@{=}[d]&A \ar[r]^{h_0}\ar[d]^{a} & C \ar[r]^{h_1} & X_2 \ar[r]^{h_2} & \cdots \ar[r]^{h_{n-1}} & X_{n} \ar[r]^{h_{n}}&X_{n+1}\ar@{=}[d] \\
\Sigma^{-1}X_{n+1}\ar[r]^{\quad d_{-1}}&B \ar[r]^{d_0} & Y_1 \ar[r]^{d_1} & Y_2\ar[r]^{d_2} & \cdots \ar[r]^{d_{n-1}} & Y_n \ar[r]^{d_{n}} & X_{n+1} }$$
where $h_{-1}=(-1)^n\Sigma^{-1}h_{n+1}$ and $d_{-1}=(-1)^n\Sigma^{-1}d_{n+1}$.
By \cite[Lemma 4.1]{BT}, there are $\varphi_1,\varphi_2,\cdots,\varphi_{n}$ which
give a morphism of $(n+2)$-angles
$$\xymatrix{
\Sigma^{-1}X_{n+1}\ar[r]^{\quad h_{-1}}\ar@{=}[d]&A \ar[r]^{h_0}\ar[d]^{a} & C \ar[r]^{h_1}\ar@{-->}[d]^{\varphi_1} & X_2 \ar[r]^{h_2} \ar@{-->}[d]^{\varphi_2}& \cdots \ar[r]^{h_{n-1}} & X_{n} \ar[r]^{h_{n}}\ar@{-->}[d]^{\varphi_n}&X_{n+1}\ar@{=}[d] \\
\Sigma^{-1}X_{n+1}\ar[r]^{\quad d_{-1}}&B \ar[r]^{d_0} & Y_1 \ar[r]^{d_1} & Y_2\ar[r]^{d_2} & \cdots \ar[r]^{d_{n-1}} & Y_n \ar[r]^{d_{n}} & X_{n+1} }$$
and the sequence
$$A\xrightarrow{\left(
              \begin{smallmatrix}
               -h_0 \\a
              \end{smallmatrix}
            \right)}C\oplus B\xrightarrow{~}X_2\oplus Y_1\xrightarrow{~}
\cdots\xrightarrow{~}X_{n}\oplus Y_{n-1}\xrightarrow{~}Y_{n}\xrightarrow{~}\Sigma A$$
is an $(n+2)$-angle.
We observe that $X_2\oplus Y_1,\cdots,X_{n}\oplus Y_{n-1},Y_{n}\in\A$.
This shows that $\left(
              \begin{smallmatrix}
               -h_0 \\a
              \end{smallmatrix}
            \right)$ is an $\s_{\A}$-inflation.
Since $\left(
              \begin{smallmatrix}
               0&1\\
                -1&0
              \end{smallmatrix}
            \right)\colon C\oplus B\to B\oplus C$ is an isomorphism, it is in particular an
an $\s_{\A}$-inflation. Thus
$\left(\begin{smallmatrix}
              a \\h_0
              \end{smallmatrix}
            \right)=\left(
              \begin{smallmatrix}
               0&1\\
                -1&0
              \end{smallmatrix}
            \right)\left(
              \begin{smallmatrix}
               -h_0 \\a
              \end{smallmatrix}
            \right)$  is also an $\s_{\A}$-inflation.

Since $\left(
              \begin{smallmatrix}
               1&0\\
                -b&1
              \end{smallmatrix}
            \right)\colon B\oplus C\to B\oplus C$ is an isomorphism, it is in particular an
an $\s_{\A}$-inflation. Thus
$\left(\begin{smallmatrix}
              a \\0
              \end{smallmatrix}
            \right)=\left(
              \begin{smallmatrix}
               1&0\\
                -b&0
              \end{smallmatrix}
            \right)\left(
              \begin{smallmatrix}
               a \\h_0
              \end{smallmatrix}
            \right)$ is also an $\s_{\A}$-inflation.
By definition there exists an $(n+2)$-angle
$$A_{\bullet}:~A\xrightarrow{\binom{a}{0}}B\oplus C\xrightarrow{(f_1,\hspace{0.8mm} q_1)}A_2\xrightarrow{f_2}A_3\xrightarrow{f_3}\cdots\xrightarrow{f_{n-1}}A_n\xrightarrow{f_n}A_{n+1}\xrightarrow{f_{n+1}}\Sigma A_0$$
with $A_2,A_3,\cdots,A_{n+1}\in\A$. By Lemma \ref{y1}, we obtain $A_{\bullet}\simeq A'_{\bullet}\oplus C_{\bullet}$, where
$$A'_{\bullet}:~A\xrightarrow{~a~}B\xrightarrow{f_{11}}A'_2\xrightarrow{f_{21}}A_3\xrightarrow{f_3}\cdots\xrightarrow{f_{n-1}}A_n\xrightarrow{f_n}A_{n+1}\xrightarrow{f_{n+1}}\Sigma A_0$$
and
$$C_{\bullet}:~0\xrightarrow{}C\xrightarrow{~1~}C\xrightarrow{}0\xrightarrow{}\cdots\xrightarrow{}0\xrightarrow{}0\xrightarrow{}0$$
is two an $(n+2)$-angle.
Since $A_2\in\A$, we also have $A'_2\in\A$.
This shows that $a$ is also an $\s_{\A}$-inflation.

Similarly, we also can prove that if $b\circ a$ is an $\s_{\A}$-deflations, then so is $b$.

By Proposition \ref{prop2}, we have that $(\A,\mathscr{E}_{\A})$
is an $n$-exact category.   \qed

\begin{remark}
For the proof of Theorem \ref{cor1}, our proof method is to avoid proving that $\mathscr{E}_{\A}$ is closed under weak isomorphisms, which is very complicated and difficult.
\end{remark}

\section*{Acknowledgments}
The author would like to thank Yonggang Hu and Tiwei Zhao for the helpful discussions.

Panyue Zhou\\
College of Mathematics, Hunan Institute of Science and Technology, 414006 Yueyang, Hunan,  People's Republic of China.\\
E-mail: \textsf{panyuezhou@163.com}\\[0.3cm]

\end{document}